\documentclass{amsart}
\usepackage{amssymb,amsmath,amsfonts,epsfig,latexsym}

\newtheorem{theorem}{Theorem}

\newtheorem{proposition}{Proposition}

\newtheorem{example}{Example}

\def\Z{\ensuremath{\mathbb{Z}}}
\def\Q{\ensuremath{\mathbb{Q}}} 
\def\P{\ensuremath{\mathbb{P}}}

\def\A{\ensuremath{\mathbb{A}}}
\def\R{\ensuremath{\mathbb{R}}}
\def\G{\mathbb{G}}
\def\l{\ell}

\def\O{\ensuremath{\mathcal{O}}}

\def\conv{\mathrm{conv}}

\def\Mov{\mathrm{Mov}}
\def\OMov{\overline{\Mov}}
\def\map{\dasharrow}

\def\int{\mathrm{int}}

\def\Span{\mathrm{Span}}

\def\OAmp{\overline{\mathrm{Amp}}\,\!}

\def\<{\ensuremath{\langle}}
\def\>{\ensuremath{\rangle}}

\begin{document}
    
\title[Stable base loci, movable curves, and small modifications]{Stable base loci, movable curves, and small modifications, for toric varieties}           

    
\author{Sam Payne}
\address{University of Michigan, Department of Mathematics, 2074 East Hall, 530 Church St., Ann Arbor, MI 48109}
\email{sdpayne@umich.edu}
\thanks{Supported by a Graduate Research Fellowship from the NSF}

\begin{abstract}
We show that the dual of the cone of divisors on a complete $\Q$-factorial toric variety $X$ whose stable base loci have dimension less than $k$ is generated by curves on small modifications of $X$ that move in families sweeping out the birational transforms of $k$-dimensional subvarieties of $X$.  We give an example showing that it does not suffice to consider curves on $X$ itself.
 \end{abstract}
 
\maketitle

\section{Introduction}

Let $X$ be a smooth $n$-dimensional complex projective variety.  Boucksom, Demailly, Paun, and Peternell have recently shown that the dual of the cone of numerical classes of effective divisors on $X$ is the closure of the cone of numerical classes of movable curves \cite[Theorem 0.2]{BDPP}; see also \cite[11.4.C]{PAG}.  Their theorem is analogous to a theorem of Kleiman that says that the dual of the cone of numerical classes of ample divisors is the closure of the cone of numerical classes of effective curves \cite[1.4.23]{PAG}.  Debarre and Lazarsfeld asked whether these results might generalize in the following way: for $1 < k < n$, is the dual of the cone of numerical classes of divisors on $X$ whose stable base locus has dimension less than $k$ generated by the numerical class of a natural collection of curves, such as curves moving in a family that sweeps out a subvariety of dimension $k$?  This paper gives an affirmative answer to their question in the toric case, with a slight twist---the curves that one must consider include not only curves on $X$, but also curves on small modifications of $X$.  

Recall that a small modification $f : X \map X^\dagger$ is a birational map that is an isomorphism in codimension 1.  If $V$ is a subvariety of $X$, we say that $f$ maps $V$ birationally to $f(V)$ if $f$ is defined on an open set meeting $V$ and maps a dense open subset of $V$ isomorphically onto a dense open subset of $f(V)$.  If $D = \sum d_i D_i$ is a divisor on $X$, then we write $f(D)$ for the divisor $\sum d_i f(D_i)$ on $X^\dagger$.

\begin{theorem}\label{main}
Let $X$ be a complete $\Q$-factorial toric variety with dense torus $T$, and let $D$ be a divisor on $X$.  The following are equivalent:
\begin{enumerate}

\item The stable base locus $B(D)$ has dimension less than $k$.

\vspace{2 pt}

\item For every $T$-invariant subvariety $V \subset X$ of dimension $k$, and for every small modification $f : X \map X^\dagger$ that maps $V$ birationally to $f(V)$, with $X^\dagger$ projective and $\Q$-factorial, and for every irreducible curve $C$ on $X^\dagger$ moving in a family sweeping out $f(V)$, $(f(D) \cdot C) \geq 0.$
\end{enumerate}
\end{theorem}

Recall that the stable base locus $B(D)$ of a $\Q$-divisor $D$ is the set-theoretic intersection of the base loci of the complete linear systems $|mD|$ for all positive integers $m$ such that $mD$ is integral.    Let $N^1(X)_\Q$ be the space of $\Q$-divisors on $X$, modulo numerical equivalence. Since the stable base locus of the sum of two divisors is contained in the union of their stable base loci, the set of classes in $N^1(X)_\Q$
that are represented by divisors whose stable base locus has codimension greater than $k$ is a convex cone (for a general numerical class, the stable base locus is independent of the choice of representative \cite[2.1 and 10.3]{PAG}).  We write $\OAmp^{k}(X)$ for the closure of this cone in $N^1(X)_\R := N^1(X)_\Q \otimes \R$. This notation is suggestive of the fact that $\OAmp^{k}(X)$ is also the closure of the cone of numerical classes of divisors that are ample in codimension $k$, that is, the cone of numerical classes of divisors $D$ such that there is an open set $U$, with $\O(D)|_U$ ample, whose complement has codimension greater than $k$.  The cone $\overline{\mathrm{Eff}}\,\!^1(X)$ of pseudo-effective divisors, the closure of the cone of numerical clases of effective divisors, is filtered by the cones $\OAmp^k(X)$:
\[
    \overline{\mathrm{Eff}}\,\!^1(X) = \OAmp^0(X) \supset \cdots \supset \OAmp^{n-1}(X) = \OAmp(X).
\]

We have an analogous filtration of the cone $\overline{\mathrm{Eff}}_1(X)_\R$ of pseudo-effective curves, as follows.  Let $N_1(X)_\R$ be the space of 1-cycles in $X$ with coefficients in $\R$, modulo numerical equivalence.  The intersection pairing on $X$ induces a natural isomorphism $N_1(X)_\R \cong N^1(X)_\R^\vee$.  Let $\Mov_k(X)$ be the  the convex cone in $N_1(X)_\R$ generated by classes of irreducible curves moving in a family that sweeps out a $k$-dimensional subvariety of $X$.  So
\[
\OMov(X) = \OMov_n(X) \subset \cdots \subset \OMov_1(X) = \overline{\mathrm{Eff}}_1(X).
\]

If a curve $C$ moves in a family that sweeps out a $k$-dimensional subvariety, and if $B(D)$ has dimension smaller than $k$, then $C$ can be moved off of the support of some effective divisor linearly equivalent to a multiple of $D$, and hence $(D \cdot C)$ is nonnegative.  So there are natural inclusions 
\[
\OMov_{k}(X) \subset \OAmp^{n-k}(X)^\vee \mbox{ \ \ and \ \ } \OAmp^{n-k}(X) \subset \OMov_{k}(X)^\vee.
\]
Kleiman's theorem says that equality holds for $k = 1$, and the theorem of Boucksom, Demailly, Paun, and Peternell says that equality holds for $k = n$.  As Lazarsfeld points out, one does not expect equality for $1 < k < n$; in light of the equality for $k = n$, equality for general $k$ would say, roughly speaking, that the restriction of a big divisor to its stable base locus is not pseudo-effective, which seems too strong (a classical theorem of Fujita says that the restriction of a divisor to its stable base locus is not ample \cite[Theorem 1.19]{Fujita}).  

Another heuristic reason to believe that equality should not hold for general $k$ comes from the consideration of small modifications $f: X \map X^\dagger$, with $X^\dagger$ complete and $\Q$-factorial.  Numerical equivalence on $X$ and on $X^\dagger$ coincide (this follows, for instance, from Kleiman's characterization of numerically trivial line bundles \cite[XIII, Theorem 4.6]{SGA6}).  So the map $D \mapsto f(D)$ induces an isomorphism $N^1(X)_\R \cong N^1(X^\dagger)_\R$.  The global sections of $D$ and of $f(D)$ are canonically identified, so if $V$ is contained in $B(D)$, and if $f$ maps $V$ birationally to $f(V)$, then $f(V)$ is in the stable base locus of $f(D)$.  In particular, $\OAmp^1(X)$ is identified with $\OAmp^1(X^\dagger)$.  But there is no natural way to identify $\OMov_{n-1}(X)$ and $\OMov_{n-1}(X^\dagger)$ compatibly with the intersection pairing.

We now introduce additional cones of numerical classes of 1-cycles which take into account $k$-movable curves on small modifications of $X$.  Let $f : X \map X^\dagger$ be a small modification, with $X^\dagger$ complete and $\Q$-factorial.  We define
\[
\Mov_k(X, X^\dagger) \subset N_1(X)_\R
\]
to be the image of the convex cone generated by numerical classes of irreducible curves $C$ on $X^\dagger$ moving in a family that sweeps out the birational image of a $k$-dimensional subvariety of $X$ under the identifications
\[
N_1(X^\dagger)_\R = N^1(X^\dagger)_\R^\vee \cong N^1(X)_\R^\vee = N_1(X)_\R.
\]
There is an obvious inclusion
\[
\OMov_k(X, X^\dagger) \subset \OAmp^{n-k}(X)^\vee.
\]
With this notation, Theorem \ref{main} can be restated as follows.

\begin{theorem}\label{main'}
Let $X$ be a complete $n$-dimensional $\Q$-factorial toric variety.  Then
\[
\OAmp^{n-k}(X)^\vee \ = \sum_{f : X \map X^\dagger} \OMov_k(X, X^\dagger),
\]
where the sum is over all small modifications $f: X \map X^\dagger$ such that $X^\dagger$ is projective and $\Q$-factorial.
\end{theorem}

It seems unclear whether a statement like Theorem \ref{main'} should be true in general, because it is not clear whether a general variety $X$ will have enough $\Q$-factorial projective small modifications.  By \cite[Lemma 1.6]{HK}, this depends essentially on whether $X$ has enough divisors $D$ whose section ring
\[
R(X,D) = \bigoplus_{m \geq 0} H^0(X, \O(mD))
\]
is finitely generated.

\vspace{10 pt}

\noindent \textbf{Problem.}  Let $X$ be a projective $n$-dimensional $\Q$-factorial complex variety.  Is the inclusion
\[
\sum_{f : X \map X^\dagger} \OMov_k(X, X^\dagger) \subset \ \OAmp^{n-k}(X)^\vee
\]
an equality?

\vspace{10 pt}

The following example illustrates the necessity of considering curves on small modifications; we give an example of a threefold $X$ such that $\OMov_2(X)$ is properly contained in $\OAmp^1(X)^\vee$.

\begin{example} \label{projective bundle} \emph{Let $Y$ be the projectivized vector bundle
\[
Y = \P(\O(3) \oplus \O_{\P^2}) \stackrel{\pi}{\longrightarrow} \P^2, 
\]
and let $s$ be the section of $\pi$ whose image is $\P(\O(3))$.  Let $p_1, p_2,$ and $p_3$ be noncolinear points in $\P^2$, and let $X$ be the blow up of $Y$ at $s(p_1)$, $s(p_2)$, and $s(p_3)$.  Let $E_i$ be the exceptional divisor over $s(p_i)$.  We make the following claims, which will be justified in the next section using toric methods.
\begin{enumerate}
\item The divisors $E_1, E_2,$ and $E_3$, together with $D^+$, the strict transform of $\P(\O(3))$, and $D^-$, the strict transform of $\P(\O_{\P^2})$, give a basis for $N^1(X)_\R$.
\nolinebreak \vspace{2 pt}
\item The class $c \in N_1(X)_\R$ defined by
$(E_i \cdot c) = 1$, $(D^+ \cdot c) = 0$, and $(D^- \cdot c) = -3$,
generates an extremal ray of $\OAmp^1(X)^\vee$.
\end{enumerate} }

\emph{ Assume these claims are true.  Choose in $N^1(X)_\R$ an open cone $U$ containing $c$ sufficiently small so that any class in $U$ pairs negatively with $D^-$ and positively with the $E_i$.  Any effective representative of a class in $U$ must have a component that lies in $D^-$,  and $D^-$ is disjoint from the $E_i$, so $U$ does not contain the class of any irreducible curve on $X$.  By perturbing a supporting hyperplane cutting out the extremal ray spanned by $c$, we can find a closed half space $H$ containing $\OAmp^1(X)^\vee \smallsetminus (\OAmp^1(X)^\vee \cap U)$ but not containing $c$.  Since $\OMov_2(X)$ is generated by classes in $\OAmp^1(X)^\vee \smallsetminus (\OAmp^1(X)^\vee \cap U)$, $\OMov_2(X)$ is contained in $H$ and does not contain $c$. }

\emph{However, consider the small modification
\[
f: X \map X^\dagger
\]
given by flopping the three curves $\pi^{-1}(p_i)$.  The variety $X^\dagger$ can be realized as a projectivized vector bundle over the blow up of $\P^2$ at $p_1$, $p_2$, and $p_3$,
\[
X^\dagger = \P(L \oplus \O_{B\l_3\P^2}) \stackrel{\pi^\dagger}{\longrightarrow} B\l_3\P^2,
\]
where $L$ is the line bundle given by the strict transform of a smooth cubic passing through the $p_i$.  Let $s^\dagger$ be the section of $\pi^\dagger$ whose image is $\P(\O_{B\l_3\P^2})$.  Then $c$ is represented on $X^\dagger$ by $s^\dagger(C)$, where $C$ is the strict transform of a smooth conic in $\P^2$ passing through the $p_i$, and $s^\dagger(C)$ moves in a family that sweeps out the surface $\P(\O_{B\l_3\P^2})$ in $X^\dagger$. }
\end{example}

\section{Toric preliminaries}

For the remainder of this paper, we fix a complete simplicial fan $\Delta$ in an $n$-dimensional latticed vector space $N_\R = N \otimes_\Z \R$.  Let $X = X (\Delta)$ be the corresponding toric variety over some fixed algebraically closed field, and let $T$ be the dense torus in $X$.  Let $\rho_1, \ldots, \rho_r$ be the rays (or 1-dimensional cones) of $\Delta$.  Let $v_i \in N$ be the primitive generator of $\rho_i$, and let $D_i$ be the corresponding prime $T$-invariant divisor on $X$.  See \cite{Fulton} for details on the correspondence between fans and toric varieties.

Since $\Delta$ is complete and simplicial, $X(\Delta)$ is $\Q$-factorial, and numerical and rational equivalence coincide on $X$.  We have a short exact sequence
\[
  0 \rightarrow N_1(X)_\R \rightarrow \R^r \rightarrow N_\R \rightarrow 0,
\]  
where the map from $\R^r$ takes $(a_1, \ldots, a_r)$ to $a_1 v_1 + \cdots + a_r v_r$.  The dual short exact sequence is
\[
0 \leftarrow N^1(X)_\R \leftarrow \R^r \leftarrow M_\R \leftarrow 0,
\]
where $M$ is the character lattice of $T$,  and the map from $\R^r$ takes $(d_1, \ldots, d_r)$ to the numerical class of the divisor $d_1 D_1 + \cdots + d_r D_r$.
  
If $D = d_1 D_1 + \cdots + d_r D_r$, then the polytope $P_D \subset M_\R$ is defined by
\[
P_D = \{ u \in M_\Q : \<u, v_i \> \geq -d_i \mbox{ for all } 1 \leq i \leq r \}.
\]
If $D$ is integral, the characters $\chi^u$, for all $u \in P_D \cap M$, form a basis for $H^0(X, \O(D))$.  The vanishing locus of $\chi^u \in H^0(X, \O(D))$ is supported exactly on the union of the divisors $D_i$ such that $\<u, v_i \> > -d_i$.  It follows that the base locus, and hence the stable base locus, of any divisor on $X$ is $T$-invariant.  See \cite{HKP2} for more details.

If $V(\tau)$ is the $T$-invariant subvariety corresponding to a cone $\tau \in \Delta$, then $V(\tau)$ is not contained in the stable base locus of $D$ if and only if the class of $D$ in $N^1(X)_\R$ is in the cone $\Gamma_\tau$ spanned by the classes of $T$-invariant divisors not containing $V(\tau)$.  In other words,
\[
\Gamma_\tau = \< D_j : \rho_j \not \in \tau \>.
\]
Since base loci are $T$-invariant,  and the codimension of $V(\tau)$ is the dimension of $\tau$,
\[
\OAmp^k(X) = \bigcap_{\dim \tau = k} \Gamma_\tau.
\]
Since each $\Gamma_\tau$ is a rational polyhedral cone, it follows that $\OAmp^k(X)$ is also a rational polyhedral cone.  Furthermore, the dual of $\OAmp^k(X)$ is given by
\[
\OAmp^k(X)^\vee = \sum_{\dim \tau = k} \Gamma_\tau^\vee.
\]
From the description of $\Gamma_\tau$ and the exact sequences above, we have
\[
\Gamma_\tau^\vee = \{ (a_1, \ldots, a_r) \in \R^r : a_1 v_1 + \cdots + a_r v_r = 0 \mbox{ and } a_i \geq 0 \mbox{ for } v_i \not \in \tau \}.
\]
The cones $\{ \Gamma_\tau \}$, for maximal cones $\tau$, appear prominently in the ``bunches of cones" in \cite{BH}.

We now prove the claims made in Example \ref{projective bundle}.

\begin{example} \emph{
Suppose $N = \Z^3$, $r = 8$, and
\[
\begin{array}{llll}
    v_{1} = (1,1,-1), & v_{2} = (-1,0,-1), & v_{3} = (0,-1,-1), & v_{4} =
    (1,0,-1), \\ v_{5} = (0,1,-1), & v_{6} = (-1,-1,-1), & v_{7} = (0, 0,
    -1), & v_{8} = (0,0,1),
\end{array}
\]
with $\Delta$ being the fan whose maximal cones are
\[
\begin{array}{llll}
    \< v_{1}, v_{4}, v_{8} \>, & \< v_{1}, v_{5}, v_{8} \>, & 
    \< v_{2}, v_{5}, v_{8} \>, & \< v_{2}, v_{6}, v_{8} \>, \\
    \< v_{3}, v_{6}, v_{8} \>, & \< v_{3}, v_{4}, v_{8} \>, &
    \< v_{1}, v_{4}, v_{5} \>, & \< v_{2}, v_{5}, v_{6} \>, \\
    \< v_3, v_4, v_6 \>, & \< v_4, v_5, v_7 \>, &
    \< v_5, v_6, v_7 \>, & \< v_4, v_6, v_7 \>. 
\end{array}
\]
Then $X = X(\Delta)$ is the variety considered in Example \ref{projective bundle}, $D_i = E_i$ for $i = 1, 2,$ and $3$, $D_7 = D^-$ is the strict transform of $\P(\O_{\P^2})$, and $D_8 = D^+$ is the strict transform of $\P(\O(3))$.  See \cite[pp.58--59]{Oda} for the construction of the fan corresponding to a projectivized split vector bundle on a toric variety.  The following diagram illustrates the fan $\Delta^\dagger$ corresponding to $X^\dagger$, as well as the fan $\Delta$; the intersection of each fan with the hyperplane $\{ (x,y,z) \in \R^3 : z = -1 \}$ is shown.  See \cite{Reid} for the changes in the fan corresponding to a flop. }

\vspace{10 pt}

\includegraphics{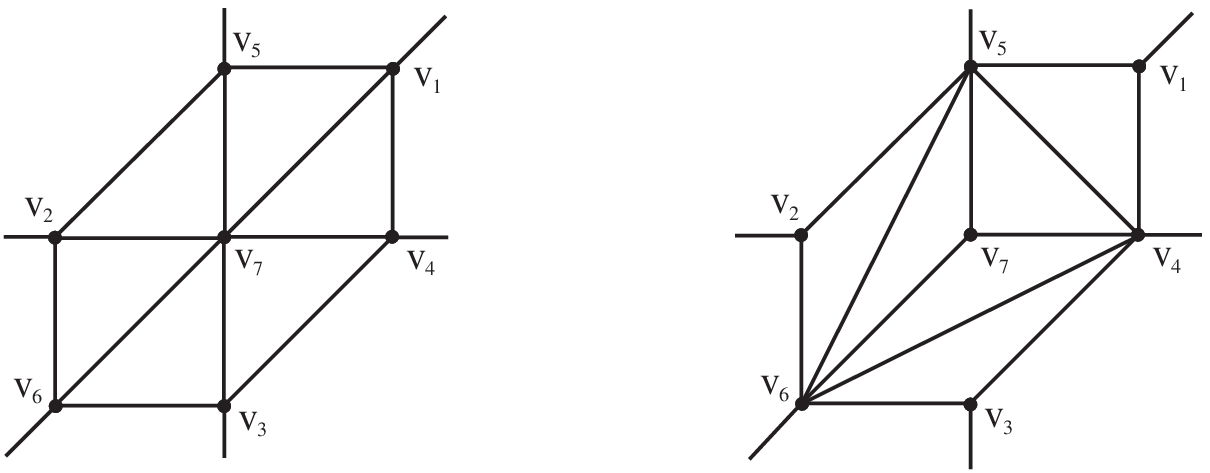}

\vspace{5 pt}
\begin{center}
$\Delta^\dagger \cap (z = -1)$ \hspace{150 pt} $\Delta \cap (z = -1)$
\end{center}
\vspace{10 pt}

\emph{ We claimed that the classes of $D_1, D_2, D_3, D_7$, and $D_8$ give a basis for $N^1(X)_\R$.  This now follows immediately from the short exact sequence
\[
0 \rightarrow M_\R \rightarrow \R^8 \rightarrow N^1(X)_\R \rightarrow 0,
\]
and the fact that $\{ v_4, v_5, v_6\}$ is a basis for $N_\R = M_\R^\vee$.  It remains to show that the class $c \in N_1(X)_\R$ given by $( D_1 \cdot c) = (D_2 \cdot c) = ( D_3 \cdot c) = 1$, $(D_7 \cdot c) = -3$, and $(D_8 \cdot c) = 0$ spans an extremal in $\OAmp^1(X)^\vee$.  Recall that
\[
\OAmp^1(X)^\vee = \sum_{i=1}^8 \Gamma_{\rho_i}^\vee.
\]
For $i \neq 7$, $D_7$ is in $\Gamma_{\rho_i}$, and hence $c$ is not in $\Gamma_{\rho_i}^\vee$.  Therefore, it suffices to show that $c$ spans an extremal ray in $\Gamma_{\rho_7}^\vee$, which is clear: $c$ spans the inward normal to the face of $\Gamma_{\rho_7}$ spanned by the classes of $D_4, D_5, D_6$, and $D_8$. }
\end{example}

\section{Constructing curves on toric varieties}

By the exact sequence
\[
0 \rightarrow N_1(X)_\R \rightarrow \R^r \rightarrow N_\R \rightarrow 0,
\]
a numerical class $c \in N_1(X)_\R$ is represented uniquely by an $r$-tuple of real numbers $(a_1, \ldots, a_r)$, such that $(D_i \cdot c) = a_i$, and these numbers satisfy the relation
\[
a_1v_1 + \cdots + a_r v_r = 0.
\] 
Conversely, any such $r$-tuple $(a_i)$ corresponds to some numerical class in $N_1(X)_\R$.

Given a cone $\tau \in \Delta$, we say that a ray $\rho_i$ is adjacent to $\tau$ if there is a cone $\sigma \in \Delta$ containing both $\tau$ and $\rho_i$.

\begin{proposition}\label{curves move}
Let $C$ be an irreducible curve in $X$, with $(D_i \cdot C) = a_i$.  Then there is a cone $\tau \in \Delta$ such that 
\begin{enumerate}
\item $C$ moves in a family sweeping out $V(\tau)$,
\item $a_i = 0$ unless $\rho_i$ is adjacent to $\tau$,
\item $a_i \geq 0$ unless $\rho_i$ is in $\tau$.
\end{enumerate}
\end{proposition}

\begin{proof}
Let $\tau$ be the unique cone such that the orbit $O_\tau$ contains an open dense subset of $C$.  The action of $T$ moves $C$ in a family sweeping out $V(\tau)$.  Furthermore, if $\rho_i$ is not adjacent to $\tau$, then $D_i$ is disjoint from $V(\tau)$, so $D_i \cdot C = 0$.  Similarly, if $\rho_i$ is not in $\tau$, then $D_i$ does not contain $O_\tau$, and hence does not contain $C$, so $(D_i \cdot C) \geq 0$.
\end{proof}

Proposition \ref{curves move} restricts the numerical classes of 1-cycles that can be represented by a positive scalar multiple of an irreducible curve.  We now show, by explicitly constructing curves with given numerical properties, that these are the only nontrivial restrictions.

\begin{proposition} \label{constructing curves}
Let $\tau$ be a cone in $\Delta$, and let $a_1, \ldots, a_r$ be integers such that
\begin{enumerate}
\item $a_1 v_1 + \cdots + a_r v_r = 0$,
\item $a_i = 0$ unless $\rho_i$ is adjacent to $\tau$,
\item $a_i \geq 0$ unless $\rho_i$ is in $\tau$.
\end{enumerate}
Then there is an irreducible curve $C$ that moves in a family sweeping out $V(\tau)$ such that $(D_i \cdot C) = a_i$ for all $1 \leq i \leq r$.
\end{proposition}

\begin{proof}
We first consider the case $\tau = 0$.  So we are given nonnegative integers $a_1, \ldots, a_r$ such that $a_1 v_1 + \cdots + a_r v_r = 0$, and must construct a curve $C$ sweeping out $X$ such that $(D_i \cdot C) = a_i$ for all $1 \leq i \leq r$.  Choose distinct elements $\lambda_1, \ldots, \lambda_r$ in the base field, let $\varphi_i : \G_m \rightarrow T$ be the one-parameter subgroup corresponding to $v_i$, and consider the rational map $\phi: \A^1 \map T$ given by
\[
\phi(z) =  \prod_{i=1}^r \varphi_i(z - \lambda_i)^{a_i},
\]
where the product is given by the group law on $T$.  Since $X$ is complete, $\phi$ extends to a regular morphism $\P^1 \rightarrow X$.  It is straightforward to check, using local coordinates, that if $a_i = 0$ then $\phi(\P^1)$ is disjoint from $D_i$, and if $a_i$ is positive then 
\[
D_i \cap \phi(\P^1) = \phi(\lambda_i).
\]
Furthermore, in the latter case, $\phi(\lambda_i)$ is a point in the dense orbit $O_{\rho_i} \subset D_i$, along which $D_i$ is locally principal, and the local intersection multiplicity is
\[
(D_i \cdot \phi(\P^1))_{\phi(\lambda_i)} = a_i.
\]
The action of $T$ moves $\phi(\P^1)$ in a family sweeping out $X$.

We now consider the case of a nonzero cone $\tau$ in $\Delta$.  Let $N_\tau$ be the quotient lattice 
\[
N_\tau = N / ( N \cap \Span\, \tau),
\]
and let $\Delta_\tau$ be the fan in $(N_\tau)_\R$ whose cones are the projections of cones of $\Delta$ containing $\tau$, so $V(\tau)$ is the toric variety corresponding to $\Delta_\tau$ \cite[3.1]{Fulton}.  The rays of $\Delta_\tau$ are the projections of the rays of $\Delta$ adjacent to $\tau$, but not contained in $\tau$; after renumbering, we may assume that these rays are $\rho_1, \ldots, \rho_s$.  Let $\overline{\rho}_i$ be the image of $\rho_i$ under the projection $\pi: N_\R \rightarrow (N_\tau)_\R$, let $w_i$ be the primitive generator of $\overline{\rho}_i$ in $N_\tau$, let $m_i$ be the positive integer such that $ \pi(v_i) = m_i w_i$, and let $\overline{D}_i$ be the divisor in $V(\tau)$ corresponding to $\overline{\rho}_i$, for $1 \leq i \leq s$.  By hypothesis,
\[
a_1 m_1 w_1 + \cdots + a_r m_r w_r  = 0
\]
in $N_\tau$.  Hence there is a curve $C$ sweeping out $V(\tau)$ such that $(\overline{D}_i \cdot C) = a_i m_i$ for $1 \leq i \leq s$.  We claim that $C$ is the required curve. Let $\iota : V(\tau) \hookrightarrow X$ be the natural inclusion, so $(D_i \cdot C) = (\iota^*D_i \cdot C)$.  By the basic properties of toric intersection theory, as developed in \cite[5.1]{Fulton}, $\iota^*D_i = \overline{D}_i / m_i$, so
\[
(D_i \cdot C) = (\overline{D}_i \cdot C) / m_i = a_i,
\]
for $1 \leq i \leq s$.  For $i > s$, if $\rho_i$ is not adjacent to $\tau$, then $D_i$ is disjoint from $V(\tau)$ and hence $D_i \cdot C = 0$.  And if $\rho_i$ is in $\tau$, then $(D_i \cdot C)$ is uniquely determined by the condition $\sum (D_i \cdot C) v_i = 0$.  Therefore $(D_i \cdot C) = a_i$, for all $i$, as required.
\end{proof}

\section{Constructing small modifications of toric varieties}

In order to prove the conjecture in the toric case, we need to know that we have enough small modifications to work with.  Although it would suffice to construct enough small modifications torically, it is perhaps helpful to observe that every small modification $f: X \map X^\dagger$, with $X^\dagger$ projective, is toric, in the sense that $X^\dagger$ has the structure of a toric variety with dense torus $T$ such that $f$ is $T$-equivariant.  This fact follows from \cite[Proposition 1.11]{HK}; it can also be seen directly by choosing a very ample divisor $A^\dagger$ on $X^\dagger$.  Then $A^\dagger$ is the birational transform of some divisor $A$ on $X$, which is linearly equivalent to a $T$-invariant divisor $A'$.  The natural map
\[
f': X \map \P(H^0(X, \O(A')))
\]
is $T$-equivariant, and agrees with the map to $X^\dagger$, embedded by $A^\dagger$, up to a projective linear change of coordinates.  

Roughly speaking, the following proposition says that, given a cone $\tau$ and a collection of rays ``close to $\tau$" and ``surrounding $\tau$" in $\Delta$, there is a complete simplicial fan $\Delta^\dagger$ containing $\tau$, whose rays are exactly the rays of $\Delta$, such that each of the rays in the given collection is adjacent to $\tau$ in $\Delta^\dagger$.  Such a fan corresponds to a complete $\Q$-factorial small modification of $X$.

\begin{proposition} \label{small modifications}
Let $v_1, \ldots, v_k$ span a cone $\tau$ in $\Delta$.  Assume there exists $s > k+1$ such that 
\begin{enumerate}
\item $\{ v_1, \ldots, v_{s-1} \}$ is linearly independent,
\item there is a linear relation
\[
a_1 v_1 + \cdots + a_k v_k = a_{k+1} v_{k+1} + \cdots + a_s v_s,
\]
with $a_i > 0$ for all $i$,
\item $v_j$ is not contained in the convex cone $\< v_1, \ldots, v_s \>$, for $j > s$.
\end{enumerate}
Then there is a complete simplicial fan $\Delta^\dagger$ containing $\tau$, whose rays are exactly the rays of $\Delta$, such that $\rho_{i}$ is adjacent to $\tau$ in $\Delta^\dagger$ for $k+1 \leq i \leq s$.  Furthermore, $\Delta^\dagger$ can be chosen such that $X(\Delta^\dagger)$ is projective.
\end{proposition}

\begin{proof}
Choose large positive numbers $p \gg q \gg 0$, and small positive numbers $\epsilon_{s+1}, \ldots, \epsilon_r$.  Let $Q$ be the polytope in $N_\R$
\[
Q = \conv \{ p v_1, \ldots, pv_k, qv_{k+1}, \ldots, qv_s, \epsilon_{s+1} v_{s+1}, \ldots, \epsilon_r v_r \},
\]
and let $\Delta_Q$ be the fan whose nonzero cones are the cones over the faces of $Q$.  It is straightforward to check that
\[
\conv \{ pv_1, \ldots, pv_k, qv_{k+1}, \ldots, \widehat{qv_i}, \ldots, qv_s \}
\]
is a face of $Q$, and so $\< v_1, \ldots, \widehat{v_i}, \ldots, v_s \>$ is a cone in $\Delta_Q$, for $k+1 \leq i \leq s$.

By construction, the rays of $\Delta_Q$ are a subset of the rays of $\Delta$, and, provided that the $\epsilon_j$ are sufficiently general, $\Delta_Q$ is simplicial.  Let $\Delta^\dagger$ be constructed by successive star subdivisions of $\Delta_Q$ with respect to the rays of $\Delta$ that are not in $\Delta_Q$.  So $\Delta^\dagger$ is also simplicial, and since, by hypothesis, none of the rays that are added lie in the cone spanned by $\{ v_1, \ldots, v_s \}$, the cones $\< v_1, \ldots, \widehat{v_i}, \ldots, v_s \>$ in $\Delta_Q$ remain unchanged in $\Delta^\dagger$.  In particular, $\tau$ is a cone in $\Delta^\dagger$, and $\rho_i$ is adjacent to $\tau$ in $\Delta^\dagger$ for $k+1 \leq i \leq s$.

It remains to check that $X(\Delta^\dagger)$ is projective.  Since $Q$ is a convex polytope, $X(\Delta_Q)$ is projective, and since $\Delta^\dagger$ is constructed from $\Delta_Q$ by a sequence of star subdivisions, $X(\Delta^\dagger)$ is constructed from $X(\Delta_Q)$, which is projective, by a sequence of blow ups.  So $X(\Delta^\dagger)$ is also projective.
\end{proof}

\section{Proof of Theorems \ref{main} and \ref{main'}}

\begin{proof} Let $c \in N_1(X)_\R$ be a numerical class spanning an extremal ray of $\OAmp^\l(X)^\vee$.  Since
\[
\OAmp^\l(X)^\vee = \sum_{\dim \sigma = \l} \Gamma_\sigma^\vee,
\]
$c$ spans an extremal ray of $\Gamma_\sigma^\vee$ for some $\l$-dimensional cone $\sigma \in \Delta$. We will show that there is a face $\tau \prec \sigma$, a small modification $f: X \map X^\dagger$ that is birational on $V = V(\tau)$, and an irreducible curve $C$ sweeping out the birational transform $f(V)$ on $X^\dagger$, such that $C$ represents the numerical class $c$.

Since $c$ is in $\Gamma_\sigma^\vee$, $c$ is given by $(D_i \cdot c) = a_i$, with $a_i \geq 0$ for $v_i \not \in \sigma$.  After renumbering, we may assume that
\[
(D_i \cdot c)  \left\{ \begin{array}{ll} < 0 & \mbox{ for } 1 \leq i \leq k. \\ > 0 & \mbox{ for } k+1 \leq i \leq s. \\ = 0 & \mbox{ for } i > s. \end{array} \right. 
\]
Let $\tau$ be the face of $\sigma$ spanned by $\{ v_1, \ldots, v_k \}$.  Now we have
\[
-a_1v_1- \cdots - a_kv_k = a_{k+1} v_{k+1} + \cdots + a_s v_s.
\]
If $k = 0$ then, by Proposition \ref{constructing curves}, there is a curve $C$ moving in a family that sweeps out $X$ such that $(D \cdot C)  = a_i$ for all $i$.  So we may assume $k \geq 1$.  Since $v_j$ is not in $\tau$ for $j > k$, we must also have $s > k + 1$.  By Propositions \ref{constructing curves} and \ref{small modifications}, it will therefore suffice to show that $\{ v_1,  \ldots, v_{s-1} \}$ is linearly independent, and that $v_j$ is not contained in the convex cone $\< v_1, \ldots, v_s \>$ for $j > s$.

Suppose there is a linear relation $b_1 v_1 + \cdots + b_{s-1} v_{s-1} = 0$. Then we have a class $b \in N_1(X)_\R$ given by
\[
(D_i \cdot b) = \left\{ \begin{array}{ll} b_i & \mbox{ for } 1 \leq i < s. \\ 0 & \mbox{ for } i \geq s. \ \end{array} \right.
\]
For small $\epsilon$, the classes $c + \epsilon b$ and $c - \epsilon b$ lie in $\Gamma_\sigma^\vee$, but not in the ray spanned by $c$, and $2c$ can be written as
\[
2c  = (c - \epsilon b) + ( c + \epsilon b),
\] 
contradicting the assumption that $c$ spans an extremal ray of $\Gamma_\sigma^\vee$.

Similarly, if $v_j$ is in the cone spanned by $\{ v_1, \ldots, v_s \}$ for some $j >s$, then we have a linear relation $v_j = b_1 v_1 + \cdots + b_s v_s$, with all of the $b_i \geq 0$.  So there is a class $b \in N_1(X)_\R$ given by 
\[
(D_i \cdot b) = \left\{ \begin{array}{ll} b_i & \mbox{ for } 1 \leq i \leq s. \\ -1 & \mbox { for } i = j. \\ 0 & \mbox{ for } i > s \mbox{ and } i \neq j. \end{array} \right.
\]
Now $b$ is contained in $\Gamma_j^\vee \subset \OAmp^{1}(X)^\vee$, which is contained in $\OAmp^\l(X)^\vee$, since $\l \geq k $ and we have assumed $k \geq 1$.  For small positive $\epsilon$, $c - \epsilon b$ is also contained in $\Gamma_\tau^\vee \subset \OAmp^\l(X)^\vee$.  Then $c$ can be written as
\[
c = (c- \epsilon b) + \epsilon b,
\]
contradicting the assumption that $c$ spans an extremal ray.
\end{proof}

\vspace{5 pt}

\noindent \textbf{Acknowledgments.}  I thank M. Hering, A. K\"uronya, and M. Musta\c t\v a for helpful conversations related to this work, and O. Debarre for useful comments on an earlier draft of this paper.  I am especially grateful to R. Lazarsfeld for suggesting this question, and for valuable advice.

\end{document}